\input amssym.def
\input amssym

\def\CC {{\Bbb C}}

\def\FF {{\Bbb F}}

\def\HH {{\Bbb H}}

\def\NN {{\Bbb N}}

\def\QQ {{\Bbb Q}}

\def\ZZ {{\Bbb Z}}

\def\Sum{\sum\limits}

\def\cross {\times}
\def\part#1#2 {{\partial {#1}/\partial {#2}}}

\def\imic{\cong}
\def\mod{\hbox{ mod }}

\def\PSL {\mathop{\rm PSL}\nolimits}
\def\SL {\mathop{\rm SL}\nolimits}
\def\Sp {\mathop{\rm Sp}\nolimits}

\def\Pic {\mathop{\rm Pic}\nolimits}

\def\Supp {\mathop{\rm Supp}\nolimits}


\def\Ker {\mathop{\rm Ker}\nolimits}


\def\diag {\mathop{\rm diag}\nolimits}

\def\Im {\mathop{\rm Im}\nolimits}

\def\hcf {\mathop{\rm gcd}\nolimits}
\def\cA {{\cal A}}

\def\cH {{\cal H}}

\def\cO {{\cal O}}
\def\cP {{\cal P}}


\def\bdv {{\bf v}}

\font\gothic=eufm10

\def\gothM {{\hbox{\gothic M}}}


\def\rationalmap{\mathrel{{\hbox{\kern2pt\vrule height2.45pt depth-2.15pt
 width2pt}\kern1pt {\vrule height2.45pt depth-2.15pt width2pt}
  \kern1pt{\vrule height2.45pt depth-2.15pt width1.8pt}\kern-1.8pt
   {\raise1.25pt\hbox{$\scriptscriptstyle\succ$}}\kern1pt}}}
\def\deep #1 {_{\lower5pt\hbox{$#1$}}}
\def\remark{\noindent{\sl Remark. }}
\def\pf{\noindent{\sl Proof: }}

\font\head=cmb10 scaled\magstep3
\font\ithead=cmti10 scaled\magstep2
\font\smallroman=cmr8
\font\smallit=cmti8
\def\qed{\vrule width5pt height5pt depth0pt\par\smallskip}
\def\surj{\to\kern-8pt\to}
\outer\def\startsection#1\par{\vskip0pt
 plus.3\vsize\penalty-100\vskip0pt
  plus-.3\vsize\bigskip\vskip\parskip\message{#1}
   \leftline{\bf#1}\nobreak\smallskip\noindent}
\outer\def\thm #1 #2 \par{\medbreak
  \noindent{\bf Theorem~#1.\enspace}{\sl#2}\par
   \ifdim\lastskip<\medskipamount \removelastskip\penalty55\medskip\fi}
\outer\def\prop #1 #2\par{\medbreak
  \noindent{\bf Proposition~#1.\enspace}{\sl#2}\par
   \ifdim\lastskip<\medskipamount \removelastskip\penalty55\medskip\fi}
\outer\def\lemma #1 #2\par{\medbreak
  \noindent{\bf Lemma~#1.\enspace}{\sl#2}\par
   \ifdim\lastskip<\medskipamount \removelastskip\penalty55\medskip\fi}
\outer\def\corollary #1 #2\par{\medbreak
  \noindent{\bf Corollary~#1.\enspace}{\sl#2}\par
   \ifdim\lastskip<\medskipamount \removelastskip\penalty55\medskip\fi}
\def\Atbil{{\cA_t^{\rm bil}}}
\def\Gtbil{{\Gamma_t^{\rm bil}}}
\def\Asixbil{{\cA_6^{\rm bil}}}
\def\Gsixbil{{\Gamma_6^{\rm bil}}}

\def\Gtnat{{\Gamma_t^\natural}}
\def\Asixnat{{\cA_6^\natural}}
\def\Gsixnat{{\Gamma_6^\natural}}
\def\Gsixnatt{{\tilde\Gamma_6^\natural}}
\def\Asixlev{{\cA_6^{\rm lev}}}
\def\Gsixlev{{\Gamma_6^{\rm lev}}}

\def\Lift{{\mathop{\rm Lift}\nolimits}}
\def\Elift{{\mathop{\rm Exp{-}Lift}\nolimits}}
\def\id{{\mathop{\rm id}\nolimits}}
\def\Zm{{\pmatrix{\tau_1&\tau_2\cr \tau_2& \tau_3}}}
\def\Div{{\mathop{\rm Div}\nolimits}}
\def\px{{\pi_6^+}}
\def\pb{{\pi_6^{\rm bil}}}
\def\pl{{\pi_6^{\rm lev}}}

\def\Eins{{\bf I}}
\baselineskip=16pt
\raggedbottom
%
%
\centerline{\head The Moduli Space  of Bilevel-6 Abelian Surfaces}
\medskip

\centerline{\ithead G.K. Sankaran \& J.G Spandaw}
\bigskip

The moduli space $\Atbil$ of $(1,t)$-polarised abelian surfaces with a
weak bilevel structure was introduced by S.~Mukai in~[Mu]. Mukai
showed that $\Atbil$ is rational for $t=2,3,4,5$. More generally, we
may ask for birational invariants, such as Kodaira dimension, of a
smooth model of a compactification of $\Atbil$: since the choice of
model does not affect birational invariants, we refer to the Kodaira
dimension, etc., of~$\Atbil$.

From the description of $\Atbil$ as a Siegel modular $3$-fold
$\Gtbil\backslash\HH_2$ and the fact that $\Gtbil\subset\Sp(4,\ZZ)$
it follows, by a result of L.~Borisov~[Bo], that $\kappa(\Atbil)=3$
for all sufficiently large~$t$. For an effective result in this
direction see~[MS]. In this note we shall prove an intermediate result
for the case~$t=6$.

\thm A The moduli space $\Asixbil$ has geometric genus
$p_g(\Asixbil)\ge 3$ and Kodaira dimension $\kappa(\Asixbil)\ge 1$.

The case $t=6$ attracts attention for two reasons: it is the first
case not covered by the results of~[Mu]; and the image of the Humbert surface
$\cH_1(1)$ in $\Atbil$, which in the cases $2\leq t\leq 5$ is a quadric
and plays an important role both in~[Mu] and below, becomes an abelian
surface (at least birationally) because the modular curve~$X(6)$ has
genus~$1$.

The method we use is that of Gritsenko, who proved a similar result
for the moduli spaces of $(1,t)$-polarised abelian surfaces with
canonical level structure for certain values of~$t$:
see~[Gr], especially Corollary~2. We use some of the weight~$3$
modular forms constructed by Gritsenko and Nikulin as lifts of Jacobi
forms in~[GN] to produce canonical forms having effective, nonzero,
divisors on a suitable projective model~$X_6$ of~$\Asixbil$. A similar
method was used by Gritsenko and Hulek in~[GH2] to give a new proof that
the Barth--Nieto threefold is Calabi-Yau.

We also derive some information about divisors in $X_6$ and linear
relations among them.

\noindent{\smallit Acknowledgements:\/ }{\smallroman We are grateful
to the DAAD and the British Council for financial assistance under ARC
Project 313-ARC-XIII-99/45.}

\startsection 1. Compactification

According to~[Mu], $\Atbil$ is isomorphic to the quotient
$\Gtbil\backslash\HH_2$, where $\HH_2$ is the Siegel upper half-plane
$\{Z\in M_{2\times 2}(\CC)\mid Z={}^\top\! Z,\ \Im Z>0\}$ and
$\Gtbil=\Gtnat\cup\zeta\Gtnat\subset\Sp(4,\ZZ)$ acts on~$\HH_2$ by
fractional linear transformations. Here $\zeta=\diag(-1,1,-1,1)$ and,
writing $\Eins_n$ for the $n\times n$ identity matrix,
$$
\Gtnat=\left\{\gamma\in\Sp(4,\ZZ)\;\left|\; \gamma-\Eins_4\in
\pmatrix{t\ZZ &  \ZZ & t\ZZ & t\ZZ  \cr
         t\ZZ & t\ZZ & t\ZZ & t^2\ZZ\cr
         t\ZZ &  \ZZ & t\ZZ & t\ZZ  \cr
          \ZZ &  \ZZ &  \ZZ & t\ZZ\cr}
\right.\right\}.
$$
We define $H(\ZZ)$ to be the Heisenberg group $\ZZ\rtimes\ZZ^2$ embedded in
$\Sp(4,\ZZ)$ as
$$
H(\ZZ)=\left\{\left.[m,n;k]=\pmatrix{1 & m & 0 & 0  \cr
         0& 1 & 0 & 0\cr
         0 & n & 1 & 0 \cr
          n & k & -m & 1\cr}\;\right|\; m,n,k \in
         \ZZ\right\}.
$$
\goodbreak

\lemma 1.1 $\Gsixnat$ is neat; that is, if $\lambda$ is an eigenvalue of
some $\gamma\in\Gsixnat$ which is a root of unity, then~$\lambda=1$.
Any torsion element of $\Gsixbil$ has order~$2$ and fixes a
divisor in~$\HH_2$.

\pf Suppose that $\gamma\in\Gsixnat$: then the characteristic polynomial
of~$\gamma$ is congruent to $(1-x)^4\mod 6$. If some
$\gamma\in\Gsixnat$ has an eigenvalue $\lambda$ which is a nontrivial
root of unity, then we may assume that~$\lambda$ is a primitive $p$th
root of unity for some prime~$p$. The minimum polynomial
$m_\lambda(x)$ of $\lambda$ over~$\ZZ$ divides the characteristic
polynomial of~$\gamma$; so $p=2,3$ or~$5$, since $\deg m_\lambda=p-1$.
But then $m_\lambda(x)=1+x$, $1+x+x^2$ or $1+x+x^2+x^3+x^4$. The
second of these does not divide $(1-x)^4$ in $\FF_2[x]$ and the other
two do not divide $(1-x)^4$ in~$\FF_3[x]$.

So any torsion element of $\Gsixbil$ is of the form
$\gamma=\zeta\gamma'$ for some $\gamma'\in\Gsixnat$; but then the
characteristic polynomial is
$$
\eqalign{
\det(\gamma-x\Eins_4)&= \det(\zeta\gamma'-x\zeta^2)\cr
                     &\equiv (1-x^2)(1+x^2)\mod 6.\cr
}
$$
From the classification of torsion elements of $\Sp(4,\ZZ)$ and their
characteristic polynomials~[Ue], it follows that $\gamma$ is conjugate
in $\Sp(4,\ZZ)$ to either $\zeta$ or $\zeta[0,1;0]$. Both these are
elements of $\Gsixbil$ of order~$2$; their fixed loci in~$\HH_2$ are
the divisors $\{\tau_2=0\}$ and $\{2\tau_2+(\tau_2^2-\tau_1\tau_3)=0\}$
respectively (Humbert surfaces of discriminants~$1$ and~$4$).~\qed

In view of Lemma~1.1, the toroidal (Voronoi, or Igusa)
compactification $(\Asixnat)^*$ of $\Asixnat=\Gsixnat\backslash\HH_2$
is smooth, cf~[SC], pp.~276--7. The action of $\zeta$ on $\Asixnat$
extends to $(\Asixnat)^*$ and the quotient~$X_6$ is a compactification
of $\Asixbil$ having only ordinary double points as
singularities. Hence $X_6$ has canonical singularities. It agrees with
the Voronoi compactification $(\Asixbil)^*$ at least in
codimension~$1$.

\startsection 2. Modular forms and canonical forms

Gritsenko and Nikulin, in~[GN], construct the weight~$3$ cusp forms
$$
\eqalign{
F_3 &= \Lift\left(\eta^5(\tau_1)\vartheta(\tau_1,2\tau_2)\right) \in
\gothM^*_3\big(\Gamma_6^+, v_\eta^8\times \id_H\big)\cr
F'_3 &= \Lift_{-1}\left(\eta^5(\tau_1)\vartheta(\tau_1,2\tau_2)\right) \in
\gothM^*_3\big(\Gamma_6^+, v_\eta^{16}\times \id_H\big)\cr
F''_3 &=
\Lift\left(\eta^3(\tau_1)\vartheta(\tau_1,\tau_2)^2\vartheta(\tau_1,2\tau_2)\right)
\in
\gothM^*_3\big(\Gamma_6^+, v_\eta^{12}\times \id_H\big)\cr
}
$$
for the extended paramodular group $\Gamma_6^+$, with character
$\chi_D$ induced from the characters $v_\eta^D\times\id_H$ of the
Jacobi group $\SL(2,\ZZ)\ltimes H(\ZZ)$. Recall (see~[GH1],~[GN]: for
compatibility with~[Mu] and other sources we work with the transposes
of the groups given in~[GN]) that $\Gamma_6^+$ is the group generated
by the paramodular group
$$
\Gamma_6=\left\{\gamma\in\Sp(4,\QQ)\;\left|\; \gamma\in
\pmatrix{\ZZ &  \ZZ & \ZZ & t\ZZ  \cr
         t\ZZ & \ZZ & t\ZZ & t\ZZ\cr
         \ZZ &  \ZZ & \ZZ & t\ZZ  \cr
          \ZZ & {{1}\over{t}} \ZZ &  \ZZ & \ZZ\cr}
\right.\right\}
$$
and the extra involution
$$
V_6=\pmatrix{0 & 1/\sqrt{6} & 0 & 0  \cr
         \sqrt{6}& 0 & 0 & 0\cr
         0 & 0 & 0 & \sqrt{6} \cr
          0 & 0 & 1/\sqrt{6} & 0\cr}.
$$

\prop 2.1 All three of $F_3$, $F'_3$ and $F''_3$ are cusp forms, without
         character, of weight~$3$ for~$\Gsixbil$.

\pf The character is induced from $v_\eta^D\times\id_H$ by the
inclusion $j:\SL(2,\ZZ)\ltimes H(\ZZ)\to \Gamma_6^+$ given by
$$
j:\left(\pmatrix{a&b\cr c&d\cr},[m,n;k]\right) \mapsto
\pmatrix{a&m&c&0\cr
         0&1&0&0\cr
         b&n&d&0\cr
         n&k&-m&1\cr}.
$$
For $\gamma\in\SL(2,\ZZ)$ we define $j_1(\gamma)=j(\gamma,[0,0; 0])$,
putting $\gamma$ in the first and third rows and columns in
$\Sp(4,\ZZ)$; and similarly $j_2(\gamma)$ puts it in the second and
fourth.

The character $v_\eta^D\times \id_H$ is trivial on $H(\ZZ)$. In the
present cases, where $D=8$, $16$ or~$12$, $v_\eta^D$ is trivial on
$\pm\Gamma(6) = \pm\Ker(\SL(2,\ZZ)\to\SL(2,\ZZ/6))$ by~[GN],
Lemma~1.2. Since $j(-\Eins_2,[0,0;0])=\zeta$, we see that
$$
\Gsixbil\cap j\big(\SL(2,\ZZ)\ltimes H(\ZZ)\big)\subseteq
j\big(\pm\Gamma(6)\ltimes H(\ZZ)\big)\subseteq\Ker\chi_D
$$ 
for $D=8,12,16$. If $D=8$ or~$16$ then, since $V_6$ and
$I=j_1\left(\pmatrix{0&1\cr -1&0\cr}\right)$ are in $\Gamma_6^+$ and
have even order and the order of $\chi_D$ is~$3$, we know that
$\chi_D(V_6)=\chi_D(I)=1$. Therefore the element
$$
J_6=IV_6IV_6=\pmatrix{0&0&-1&0\cr 0&0&0&-6\cr 1&0&0&0\cr
0&{{1}\over{6}}&0&0\cr}\in \Gamma_6^+
$$
is in $\Ker\chi_D$. If $D=12$ then $\chi_{12}(J_6)=\chi_{12}(IV_6)^2=1$
so again $J_6\in \Ker\chi_D$. Now we proceed as in~[Gr], Lemma~2.2, and
show that the group generated by $j\big(\Gamma(6)\ltimes
H(\ZZ)\big)$ and $J_6$ includes~$\Gsixnat$. To see this, we work with
the conjugate groups $\Gsixnatt=\nu_6(\Gsixnat)$
and~$\tilde\Gamma_6=\nu_6(\Gamma_6)$, where $\nu_6$ denotes
conjugation by $R_6=\diag(1,1,1,6)$. Note that $\nu_6(J_6)=
R_6J_6R_6^{-1} =\pmatrix{0&-\Eins_2\cr \Eins_2&0}$. If
$\tilde\gamma\in\Gsixnatt$ then its second row $\tilde\gamma_{2*}$ is
$(0,1,0,0)\mod 6$. Suppose first that $\tilde\gamma_{22}=1$ and put
$$
\tilde\beta = \nu_6\big(
J_6[\tilde\gamma_{21}/6,\tilde\gamma_{23}/6;
 \tilde\gamma_{24}/6]J_6^{-1}\big) 
=
 \pmatrix{1&0&0&\tilde\gamma_{23}/6\cr
\tilde\gamma_{21}&1&\tilde\gamma_{23}&\tilde\gamma_{24}\cr
0&0&1&\tilde\gamma_{23}/6\cr
0&0&0&1\cr}
$$
Now $(0,1,0,0)\tilde\beta=\tilde\gamma_{2*}$ so the second row of
$\tilde\gamma\tilde\beta^{-1}\in\Gsixnatt$ is $(0,1,0,0)$. Such a matrix
is in $\nu_6\left(j\big(\Gamma(6)\ltimes H(\ZZ)\big)\right)$.

It remains to reduce to the case $\tilde\gamma_{22}=1$. Certainly the
vector $\tilde\gamma_{2*}$ is primitive, since $\det\tilde\gamma=1$,
and since $\tilde\gamma\in\Gsixnatt$ we have $\hcf(6,
\tilde\gamma_{21}, \tilde\gamma_{23})=6$. In the proof of [FS],
Satz~2.1 it is shown that there are integers $\lambda$, $\mu$ such
that $\tilde\gamma'=\tilde\gamma \nu_6\left([\mu,0;0] J_6
[0,\lambda;0]J_6^{-1}\right)$ has
$\hcf(\tilde\gamma'_{21},\tilde\gamma'_{23})=6$, so the second row of
$\tilde\gamma'$ is $(6x_1,6x_2+1,6x_3,6x_4)$ with
$\hcf(x_1,x_3)=1$. But then the $(2,2)$-entry of
$\tilde\gamma'\nu_6([m,n;0])$ is $6(mx_1+nx_3+x_2)+1$ which is equal
to~$1$ if we choose $m$ and~$n$ suitably.~\qed

\prop 2.2 The differential forms $\tilde\omega=F_3\,
d\tau_1\wedge d\tau_2\wedge d\tau_3$, $\tilde\omega'=F'_3\,
d\tau_1\wedge d\tau_2\wedge d\tau_3$ and
$\tilde\omega''=F''_3\, d\tau_1\wedge d\tau_2\wedge
d\tau_3$ give rise to canonical forms $\omega, \omega',
\omega''\in H^0(K_{X_6})$.

\pf By Proposition~2.1, $\tilde\omega$, $\tilde\omega'$ and
$\tilde\omega''$ are all $\Gsixbil$-invariant,
so they give rise to forms $\omega$, $\omega'$, $\omega''$
on~$\Asixbil$. Since $F_3$, $F'_3$ and $F_3''$ are cusp forms,
if any of $\omega$, $\omega'$ and $\omega''$ are holomorphic on
$\Asixbil$ they extend holomorphically to the cusps
of $(\Asixbil)^*$. Since $X_6$ agrees with $(\Asixbil)^*$ in
codimension~$1$ and has canonical singularities it follows that these
forms can be thought of as $3$-forms on~$X_6$ holomorphic at
infinity. We need to check that $\omega$, $\omega'$ and $\omega''$ are
holomorphic everywhere. But this is a well-known result of
Freitag([Fr], Satz~II.2.6).~\qed

\startsection 3. Divisors in the moduli spaces.

In this section we shall describe the canonical divisors
$\Div_{X_6}(\omega)$, $\Div_{X_6}(\omega')$  and $\Div_{X_6}(\omega'')$
in~$X_6$ and give some detail about the branching locus in $X_6$ arising
from torsion in~$\Gsixbil$.

$\Gsixbil$ is a subgroup both of the paramodular group $\Gamma_6$ and of
$\Gamma^+_6$. Hence there is a finite morphism
$\sigma:\Asixbil\to\cA_6^+$. We denote the projection map
$\HH_2\to\Asixbil$ by $\pb$ and similarly $\pi_6$, $\px$, etc.

For discriminant $\Delta=1$, $4$ we put
$$
\cH_\Delta(k)= \left\{\Zm\in\HH_2 \Bigm|
{\textstyle{{1}\over{24}}}(k^2-\Delta)\tau_1 + k\tau_2 + 6\tau_3\right\} =
0 
$$
where $k\in\ZZ$ is chosen so that ${{1}\over{24}}(k^2-\Delta)\in\ZZ$.
The irreducible components of the Humbert surfaces $H_1$ and $H_4$ of
discriminants~$1$ and~$4$ in $\cA_6$ are $\pi_6(\cH_1(k))$ and
$\pi_6(\cH_4(k))$ for $0\le k< 6$: the statement of~[vdG],
Theorem~IX.2.4 and of~[GH1], Corollary~3.3, is wrong because
$\cH_\Delta(-k)$ is $\Gamma_t$-equivalent to
$\cH_\Delta(k)$. Nevertheless the irreducible components of the
Humbert surfaces of discriminants~$1$ and~$4$ in $\cA^+_6$ are as
stated in~[GN], namely $\px(\cH_1(1))$ and $\px(\cH_1(5))$ for
discriminant~$1$ and $\px(\cH_4(1))$ for discriminant~$4$.
                                                          
The calculation of the divisors uses the product expansion of the
modular forms $F_3$, $F'_3$ and $F''_3$ given in~[GN]. We have chosen
to work with the transposes of the matrices given in~[GN], so we have
to write $q=e^{2\pi i\tau_1}$, $r=e^{2\pi i\tau_2/6}$ and $s=e^{2\pi
i\tau_3/36}$ for these expansions to be correct. This is because
$^\top\Gamma_t=\diag(1,t,1,t^{-1})\Gamma_t\diag(1,t^{-1},1,t)$ (for
any $t\in \NN$), and $\diag(1,t,1,t^{-1}):(\tau_1,\tau_2,\tau_3)\to
(\tau_1,t\tau_2, t^2\tau_3)$. A similar correction is needed in~[GH2]

By [GN], equations (4.12)--(4.14) we have (correcting a minor misprint)
$$
\eqalign{
F_3&=\Elift(5\phi_{0,3}^2-4\phi_{0,2}\phi_{0,4})=\Elift(\phi_3)\cr
F'_3&=\Elift(\phi_{0,3}^2)=\Elift(\phi'_3)\cr
F''_3&=\Elift(3\phi_{0,3}^2-2\phi_{0,2}\phi_{0,4})=\Elift(\phi''_3).\cr
}
$$
By [GN], Example~2.3 and Lemma~2.5, we have
$$
\eqalign{
\phi_{0,2}&=(r^{\pm 1}+4)+q(r^{\pm 3}-8r^{\pm 2}-r^{\pm 1}+16)+O(q^2)\cr
\phi_{0,3}&=(r^{\pm 1}+2)+q(-2r^{\pm 3}-2r^{\pm 2}+2r^{\pm 1}+4)+O(q^2)\cr
\phi_{0,4}&=(r^{\pm 1}+1)+q(-r^{\pm 4}-r^{\pm 3}+r^{\pm 1}+2)+O(q^2).\cr
}
$$

\prop 3.1 The divisors in $\HH_2$ of the cusp forms are
$$
\eqalign{
\Div(F_3) &= (\px)^{-1}\left(\px \big(\cH_1(1) + 5\cH_1(5) +
\cH_4(1)\big) \right),\cr
\Div(F'_3) &= (\px)^{-1}\left(\px \big(5\cH_1(1) + \cH_1(5) +
\cH_4(1)\big) \right),\cr
\Div(F''_3) &= (\px)^{-1}\left(\px \big(3\cH_1(1) + 3\cH_1(5) +
\cH_4(1)\big) \right).\cr
}
$$

\remark This corrects the coefficients given in [GN],~Example~4.6:
for instance,
it is easy to see, by considering the effect of an element of
order~$2$ fixing an Humbert surface, that the coefficients of
$\cH_1(1)$, $\cH_1(5)$ and $\cH_4(1)$ must be odd.

\pf Write $\phi_3=\sum f(n,l)q^nr^l$, and similarly for $\phi'_3$ and
$\phi''_3$. By [GN], Theorem~2.1, the coefficient of
$\px\big(\cH_\Delta(b)\big)$ in $\cA_6^+$ is
$$
m_{\Delta,b}=\Sum_{d>0}f(d^2a,db)
$$
where $b^2-24a=\Delta$. So to calculate $m_{1,1}$ we may take $b=1$
and $a=0$, so $m_{1,1}=\Sum_{d>0}f(0,d)$. From the formulae above,
$\phi_3=(r^{\pm 2}+6)+O(q)$, so $m_{1,1}=f(0,2)=1$. Similarly we have
$\phi'_3=(r^{\pm 2}+4r^{\pm 1}+6)$ so $m'_{1,1}=5$ and
$\phi''_3=(r^{\pm 2}+2r^{\pm 1}+6)$ so $m'_{1,1}=3$.

To calculate the coefficients of $\px\big(\cH_4(1)\big)$ we note that
$\cH_4(1)$ is $\Gamma_6^+$-equivalent to $\cH_4(2)$, so we may as well
work with that and calculate~$m_{4,2}$. For this purpose we can take
$b=2$ and $a=0$; so $m_{4,2}=\Sum_{d>0}f(0,2d)=1$, and
$m'_{4,2}=m''_{4,2}=1$ also.

To calculate $m_{1,5}$ we take $b=5$ and $a=1$, so
$m_{1,5}=\Sum_{d>0}f(d^2,5d)$. The Fourier coefficient
$f(n,l)$ depends only on $24n-l^2$ and on the residue class of
$l\mod 12$ (see [GN]); that is, in our case, on $d^2$ and
on $d\mod 12$. If $d\not\equiv \pm 1\mod 6$ then $5d\equiv \pm d\mod 12$,
so $f(d^2,5d)=f(0,\pm d)$ which is zero unless $d=\pm 2$ or $d=0$.
Since we are only interested in $d>0$ the only contribution for
$d\not\equiv \pm 1\mod 6$ arises from $d=2$, when $f(4,10)=f(0,-2)=1$.
If $d\equiv\pm 5\mod 12$ then $f(d^2,5d)=f({{-d^2+1}\over{24}},\pm 1)$
which vanishes because $f(n,l)=0$ for~$n<0$. If $d\equiv \pm 1\mod 12$
then $f(d^2,5d)=f({{-d^2+25}\over{24}},\pm 5)$ which vanishes except
possibly when~$d=1$. So $m_{1,5}=1+f(1,5)$ and from the expansions
of $\phi_{0,2}$, $\phi_{0,3}$ and $\phi_{0,4}$ we calculate $f(1,5)=4$.
Similarly $m'_{1,5}=1+f'(1,5)=1$ and $m''_{1,5}=1+f''(1,5)=3$.~\qed

Brasch [Br] has studied the branch locus of $\pi^{\rm
lev}_t:\HH_2\to\cA^{\rm lev}_t$ for all~$t$: for $t\equiv 2\mod 4$ the
divisorial part has five irreducible components. They are
$\pl(\cH_{\zeta_i})$ for $0\le i\le 4$, where $\cH_{\zeta_i}\subset
\HH_2$ is the fixed locus of $\zeta_i$ and
$$
\zeta_0=\zeta,\qquad \zeta_1=\zeta^\top[-6,0;0],
\qquad\zeta_2=\pmatrix{-7&4&0&0\cr -12&7&0&0\cr 0&0&-7&-12\cr 0&0&4&7\cr},
$$
$$
\zeta_3=\zeta[1,0;0],
\qquad
\zeta_4=\pmatrix{-1&-1&0&6\cr 0&1&-6&0\cr 0&0&-1&0\cr 0&0&-1&1\cr}.
$$

These are all elements of $\Gsixbil$. Their fixed loci are
$$
\cH_{\zeta_0}=\{\tau_2=0\},\qquad
\cH_{\zeta_1}=\{6\tau_1-2\tau_2=0\},\qquad
\cH_{\zeta_2}=\{6\tau_1-7\tau_2+2\tau_3=0\},
$$
$$
\cH_{\zeta_3}=\{2\tau_2+\tau_3=0\},\qquad
\cH_{\zeta_4}=\{2\tau_2+\tau_3-6=0\},
$$
of discriminants $1,4,1,4,4$ respectively. Thus three of the components
have discriminant~$4$ and
therefore map to $\px\cH_4(1)\subset \cA^+_6$ (they correspond to
bielliptic abelian surfaces). $\cH_{\zeta_0}=\cH_1(1)$ corresponds to
product surfaces $E\cross E'$ with polarisation given by
$\cO_E(1)\boxtimes\cO_{E'}(6)$, and $\cH_{\zeta_2}$ maps to
$\px\big(\cH_1(5)\big)$, corresponding to abelian surfaces $E\cross E'$
with polarisation $\cO_E(2)\boxtimes\cO_{E'}(3)$.

\prop 3.2 The branch locus of $\pb:\HH_2\to\Asixbil$ has seven irreducible
components, each with branching of order~$2$. They are
$\pb(\cH_{\zeta_i})$ and two other components $\pb(\cH_{\zeta'_1})$,
$\pb(cH_{\zeta''_1})$, which are equivalent to $\pb(\cH_{\zeta_1})$
in~$\Asixlev$.

\pf It follows from Lemma~1.1 that the branch locus consists of
divisors only and that the branching is of order~$2$.

Write $G=\Gsixlev\vartriangleright H=\Gsixbil$ and let $G$ act on
$\Omega=G/H \imic\PSL(2,\ZZ/6)$. By [Br], Corollary~1.3, the number of
irreducible divisors in $\Asixbil$ mapping to $\pl(\cH_{\zeta_i})$,
which is equal to the number of $H$-conjugacy classes in the
$G$-conjugacy class of~$\zeta_i$, is $|G:H.C_G(\zeta_i)|$. Moreover,
for fixed~$i$, these divisors are permuted transitively by $\Omega$ so
they all have the same branching behaviour: $\pb$ is branched of
order~$2$ above each one.

$|G:H.C_G(\zeta_i)|=|G/H:C_G(\zeta_i)/H\cap C_G(\zeta_i)|$, which is
the index of the image of $C_G(\zeta_i)$ in~$\Omega$. For $i=0,1,2,3$
the centraliser $C_{\Sp(4,\QQ)}(\zeta_i)$ is described in [Br], Lemma~2.1,
and $C_G(\zeta_i)=C_{\Sp(4,\QQ)}(\zeta_i)\cap G$.

For $\zeta_0$, if $g\in\PSL(2,\ZZ/6)\imic\Omega$ and
$\tilde \gamma\in\SL(2,\ZZ)$ is some lift of~$\gamma$ then
$j(\tilde \gamma,[0,0;0])\in C_G(\zeta_0)$ so the index is~$1$.

For $\zeta_1$, if $\gamma=\pmatrix{a&b\cr c & d\cr}
\in\PSL(2,\ZZ/6)$ and $b$ is even then
$$
\pmatrix{\tilde a&0&\tilde b&3\tilde b\cr
3(\tilde a-1)&1&3\tilde b&0\cr
\tilde c&0&\tilde d&3(\tilde d-1)\cr
0&0&0&1\cr}\in
C_G(\zeta_1)
$$
for a lift $\tilde\gamma$; and this is a necessary condition for such an
element to exist since if $\beta=\beta_{ij}\in C_G(\zeta_1)$ then
$3\beta_{13}\equiv 0 \mod 6$. So $C_G(\zeta_1)/C_G(\zeta_1)\cap H\subset
\PSL(2,\ZZ/6)$ is the reduction$\mod 6$ of~${}^\top\Gamma_0(2)$, i.e. the
preimage of $\left\{\pmatrix{a&0\cr c& d\cr}\in\SL(2,\ZZ/2)\right\}$, which
is of index~$3$ because it is the stabiliser of $(1,0)$ when
$\SL(2,\ZZ/2)$ acts as the symmetric group on the nonzero vectors
in~$\FF_2^2$.

For $\zeta_2$, an element of $C_{\Sp(4,\QQ)}$ is determined (see [Br])
by two elements of $\SL(2,\QQ)$. If we take
$$
\gamma=
\pmatrix{3&4\cr 2&3\cr},\qquad
\gamma'=
\pmatrix{10&9\cr 11&10\cr}
$$
we get an element $\beta$ whose image in $\PSL(2,\ZZ/6)$ is
$\pmatrix{0&1\cr -1&0\cr}$. If we take
$$
\gamma=
\pmatrix{11&4\cr 8&3\cr},\qquad
\gamma'=
\pmatrix{7&9\cr 3&4\cr}
$$
we get an element $\beta'$ whose image in $\PSL(2,\ZZ/6)$ is
$\pmatrix{-1&-1\cr 1&0\cr}$. These two elements generate $\PSL(2,\ZZ/6)$
because their lifts generate $\SL(2,\ZZ)$. The elements
$$
\beta=\pmatrix{-18&14&25&42\cr
               -42&31&42&72\cr
               107&-70&-18&-42\cr
               -70&46&-14&31\cr}, \qquad
\beta'=\pmatrix{23&-30&25&42\cr
                24&-5&42&72\cr
                59&-34&0&-6\cr
                -34&20&-6&7\cr}
$$
both belong to $\Gsixlev$, so the index we want is~$1$.

For $\zeta_3$, as for $\zeta_0$, $j(\tilde\gamma, [0,0;0])\in C_G(\zeta_3)$
so the index is~$1$.

For $\zeta_4$, note that $\zeta_4={}^\top[0,0;6]\zeta_3(^\top[0,0;6])^{-1}$
so $C_{\Sp(4,\QQ)}(\zeta_4)=
{}^\top[0,0;6]C_{\Sp(4,\QQ)}\zeta_3(^\top[0,0;6])^{-1}$. It happens that
${}^\top[0,0;6]j(\tilde\gamma,[0,0;0])(^\top[0,0;6])^{-1}=
j(\tilde\gamma,[0,0;0])$, so again the index is~$1$.~\qed

Next we look at the boundary divisors of~$X_6$. These correspond to
$1$-dimensional subspaces of $\QQ^4$ up to the action of~$\Gsixbil$.
We may think of such a space as being given by a unique, up to sign,
primitive vector $\bdv=(v_1,v_2,v_3,v_4)\in \ZZ^4$. It is shown in
[FS] that the $\Gamma_6$-orbit of $\bdv$ is determined by
$r=\hcf(6,v_1,v_3)$, so $\cA_6$ has four corank~$1$ cusps (or boundary
divisors in the toroidal compactification). However, the cusps $r=1$
and $r=6$ are interchanged by~$V_6$, as are the cusps $r=2$ and $r=3$,
so $\cA_6^+$ has just two corank~$1$ cusps. Since $F_3$, $F'_3$
and $F''_3$ are modular forms (with character) for $\Gamma_6^+$, the
order of vanishing of any of them at a cusp of $X_6$ given by~$\bdv$
depends only on which cusp of $\cA_6^+$ it lies over, i.e. on whether
$r$ is or is not a proper divisor of~$6$.

We write $D_1$ for the divisor in $X_6$ which is the sum of all the
boundary components with $r=1$ or $r=6$, and $D_2$ for the sum of
all the components with $r=2$ or $r=3$. By an easy modification of the
argument of [FS], Satz~2.1, one can check that both $D_1$ and $D_2$ have
$24$ irreducible components, but we shall not make any use of this.

\thm 3.3 The divisors of $\omega$, $\omega'$ and $\omega''$ in $X_6$ are
$$
\eqalign{
\Div_{X_6}(\omega) &= 4\pb(\cH_{\zeta_2})+D_1+D_2,\cr
\Div_{X_6}(\omega') &= 4\pb(\cH_{\zeta_0})+3(D_1+D_2),\cr
\Div_{X_6}(\omega'') &= 2\pb(\cH_{\zeta_0})+2\pb(\cH_{\zeta_2})
    +2(D_1+D_2).\cr
}
$$

\pf If $\pb$ is branched along the irreducible
divisors $B_\alpha$ with ramification index~$e_\alpha$, then
$d\tau_1\wedge d\tau_3\wedge d\tau_3$ acquires poles of order
$e_\alpha/2$ along~$B_\alpha$. So by Proposition~3.1
$$
\eqalign{
\Div_{X_6}(\omega) &= \sigma^{-1}\px \big(\cH_1(1) +5\cH_1(5) +
\cH_4(1)\big) -{{1}\over{2}}\sum e_\alpha B_\alpha +D,\cr
\Div_{X_6}(\omega') &= \sigma^{-1}\px \big(5\cH_1(1) +\cH_1(5) +
\cH_4(1)\big)  -{{1}\over{2}}\sum e_\alpha B_\alpha +D',\cr
\Div_{X_6}(\omega'') &= \sigma^{-1}\px \big(3\cH_1(1) +3\cH_1(5) +
\cH_4(1)\big)  -{{1}\over{2}}\sum e_\alpha B_\alpha +D'',\cr
}
$$
where $D$, $D'$ and $D''$ are effective divisors supported on the
boundary~$X_6\setminus\Asixbil$. The form of the branch locus part
of the divisors follows now from Proposition~3.2 and the discriminants
of~$\cH_{\zeta_i}$.

It remains to calculate the vanishing orders of the forms at each boundary
divisor. For each form, we need only consider two boundary components,
one from $D_1$ and one from $D_2$. We use the components $D(\bdv_1$,
$D(\bdv_2)$ corresponding to $\bdv_1=(0,0,1,0)$ and $\bdv_2=(0,0,2,1)$.
The first step in constructing the toroidal compactification near
$D(\bdv_1)$ is to take a quotient by the lattice~$P'_{\bdv_1}(\Gsixbil)$
(see for instance~[GH2], pp.925--926 or for a full explanation [HKW],
Section~I.3D). As in [HKW], Proposition~I.3.98, $P'_{\bdv_1}(\Gsixbil)$
is generated by $j_1\left(\pmatrix{1&6\cr 0&1\cr}\right)$; so a local
equation for $D(\bdv_1)$ at a general point is $t_1=0$, where
$t_1=e^{2\pi i\tau_1/6}=q^{1/6}$. Using the values of $f(0,l)$ calculated
above and the Fourier expansion given in [GN], Theorem~2.1, we see that
the expansions of $F_3$, $F'_3$ and $F''_3$ begin $q^{1/3}rs^2$,
$q^{2/3}r^3s^4$ and $q^{1/2}r^2s^3$ respectively, so their orders of
vanishing along $D_1$ are $2$, $4$ and~$3$. The form
$d\tau_1\wedge d\tau_2\wedge d\tau_3$ contributes a simple pole at the
boundary so the coefficients of $D_1$ in the divisors of $\omega$, $\omega'$
and $\omega''$ are $1$, $3$ and~$2$.

We put
$$
\theta=\pmatrix{1&-1&0&0\cr -1&2&0&0\cr 0&0&2&1\cr 0&0&1&1}\in\Sp(4,\ZZ),
$$
so that $\bdv_2=\bdv_1\theta$. Then
$\cP_{\bdv_2}=\theta^{-1}\cP_{\bdv_1}\theta$ (where, as in [HKW], $\cP_\bdv$
denotes the stabiliser of $\bdv$ in $\Sp(4,\QQ)$), and from this one readily
calculates that
$$
P'_{\bdv_2}(\Gsixbil)=\left\{\left.\pmatrix{1&0&4n&2n\cr 0&1&2n&n\cr
0&0&1&0\cr 0&0&0&1}\right| n\equiv 0 \mod 36\right\}.
$$
So the cusp $D_2$ is given by $t_2=0$, where $t_2=e^{2\pi
i(\tau_1/144+\tau_2/72+\tau_3/36)}=q^{1/144}r^{1/12}s$. The number of
times this term divides the expressions for $F_3$, $F'_3$ and $F''_3$
is in fact equal to the power of $s$ that occurs, namely $2$, $4$
and~$3$ respectively; so we get the same orders of vanishing along
$D_2$ as along~$D_1$.~\qed                          

This calculation shows directly (without appealing to Freitag's
result in~[Fr]) that $\omega$, $\omega'$ and $\omega''$ are all
holomorphic.

\corollary 3.4 In $\Pic X_6$,
$D_1+D_2=2\big(\px(\cH_{\zeta_2})-\px(\cH_{\zeta_0})\big)$.

\remark Notice that $\Div_{X_6}(\omega)+\Div_{X_6}(\omega')=
2\Div_{X_6}(\omega'')$, reflecting the fact (easily seen from~[GN])
that $F_3F'_3=(F''_3)^2$.

Theorem A now follows at once from the following observation.

\prop 3.5 $\omega$, $\omega'$ and $\omega''$ are linearly independent
elements of $H^0(K_{X_6})$.

\pf Suppose that $\lambda\omega+\lambda'\omega'+\lambda''\omega''=0$.
At a general point of $\px(\cH_{\zeta_0})$, $\omega'$ and $\omega''$
vanish but $\omega$ does not. Therefore $\lambda=0$. Similarly
$\lambda'=0$, considering a general point of $\px(\cH_{\zeta_2})$. Finally,
$\omega''\neq 0$ because~$F''_3$ is not identically zero.~\qed

We want to remark that $\kappa(\Asixbil)\ge 1$ can be deduced from the
existence of $\omega'$ alone. The divisor
$\Div_{X_6}(\omega')$ is effective and
$\pb(\cH_\zeta)\subset\Supp\Div_{X_6}(\omega')$. Since $X_6$ has
canonical singularities, $K$ is effective on any smooth model of
$X_6$, and hence also on any minimal model~$X'_6$ of~$X_6$. Any
surfaces contracted by the birational map $X_6\rationalmap X'_6$ must
be birationally ruled. But $\pb(\cH_\zeta)$ is not birationally ruled:
it is isomorphic to $X(6)\times X(6)$, since $\cH_\zeta$ is isomorphic
to $\HH\times \HH$ and is preserved by the subgroup $\Gamma(6)\times
\Gamma(6)$ embedded in $\Gsixbil$ by~$(j_1,j_2)$. Thus its closure
is birationally an abelian surface, since $X(6)$ has genus~$1$. So the
canonical divisor of $X'_6$ is effective and nontrivial; so, by
abundance, some multiple of it moves and
therefore~$\kappa(\Asixbil)\geq 1$.

\beginsection References\par
\exhyphenpenalty100
\frenchspacing

\item{[SC]}A.~Ash, D.~Mumford, M.~Rapoport \& Y.~Tai, {\sl Smooth
Compactification of Locally Symmetric Varieties}, Math.\ Sci.\ Press,
Brookline, Mass., 1975.

\item{[Br]}H.-J.~Brasch, {\sl Branch points in moduli spaces of
certain abelian surfaces}, in: {\sl Abelian Varieties, Egloffstein 1993},
de Gruyter, Berlin 1995, 27--54.

\item{[Bo]}L.A.~Borisov, {\sl A finiteness theorem for subgroups of\/
$\Sp(4,\ZZ)$}, Algebraic geometry,~9. J. Math.
Sci. (New York)~{\bf 94} (1999), 1073--1099.

\item{[FS]}M.~Friedland \& G.K.~Sankaran, {\sl Das Titsgeb\"aude von
Siegelschen Modulgruppen vom Geschlecht~2}, Preprint {\tt
math.AG/0002249}, 2000, to appear in Math. Abh. Sem. Univ. Hamburg

\item{[Fr]}E.~Freitag, {\sl Siegelsche Modulfunktionen}, Grundlehren der
mathematischen Wissenschaften Bd.~{\bf 254}, Springer-Verlag, Berlin, 1983.

\item{[vdG]}G.~van~der~Geer, {\sl Hilbert Modular Surfaces}, Ergebnisse
der Mathematik und ihrer Grenzgebiete Bd.~{\bf 16}, Springer-Verlag,
Berlin-New York, 1988. 

\item{[Gr]}V.~Gritsenko, {\sl Irrationality of the moduli spaces of
polarized abelian surfaces}, in: {\sl Abelian Varieties, Egloffstein
1993}, de Gruyter, Berlin 1995, 63--84.

\item{[GH1]}V.~Gritsenko \& K.~Hulek, {\sl Minimal Siegel modular
threefolds}, Math.\ Proc.\ Cambridge Philos.\ Soc.~{\bf 123} (1998),
461--485.

\item{[GH2]}V.~Gritsenko \& K.~Hulek, {\sl The modular form of the
Barth-Nieto quintic}, Internat.\ Math.\ Res.\ Notices~{\bf 17} (1999),
915--937.

\item{[GN]}V.~Gritsenko \& V.~Nikulin, {\sl Automorphic forms and
Lorentzian Kac-Moody Algebras II}, Internat.\ J.\ Math.~{\bf 9}
(1998), 201--275.

\item{[HKW]} K.~Hulek, C.~Kahn \& S.~Weintraub, {\sl Moduli spaces of
abelian surfaces: Compactification, degenerations and theta
functions}, de Gruyter 1993.

\item{[MS]}A.~Marini \& G.K.~Sankaran, {\sl Abelian surfaces with
bilevel structure}, in preparation.

\item{[Mu]}S.~Mukai, {\sl Moduli of abelian surfaces and regular
polyhedral groups}, in: {\sl Moduli of algebraic varieties, Sapporo
1999}, to appear.

\item{[Ue]}K.~Ueno, {\sl On fibre spaces of normally
polarized abelian varieties of dimension $2$. II.  Singular fibres of
the first kind}, J.\ Fac.\ Sci.\ Univ.\ Tokyo\ Sect.\ IA\ Math.~{\bf
19} (1972), 163--199.

\bigskip
\settabs2\columns
\+G.K.~Sankaran&J.G.~Spandaw\cr
\+Department of Mathematical Sciences&Institut f\"ur Mathematik\cr
\+University of Bath&Universit\"at Hannover\cr
\+Bath BA2~7AY&Postfach 6009\cr
\+ENGLAND&D-30060 Hannover\cr
\+&GERMANY\cr
\+\cr
\+\tt gks@maths.bath.ac.uk&\tt spandaw@math.uni-hannover.de\cr
\end